\documentclass{ifacconf}

\usepackage{graphicx}      % include this line if your document contains figures
\usepackage{natbib}        % required for bibliography
\usepackage{amsmath,amssymb,amsfonts}
\usepackage{bm} % if you want bold math vectors via \bm{x}
\usepackage{float} 
\usepackage{algorithm}
\usepackage{algpseudocode}
\usepackage{pgfplots}
\usepackage{subcaption}
\usepackage{booktabs}
\usepackage{array}

\usepackage{pgfplots}
\pgfplotsset{compat=1.18}

\newcommand{\rev}[1]{\textcolor{black}{#1}}
\newcommand{\meisam}[1]{\textcolor{black}{#1}}

\makeatletter
\renewcommand\@fnsymbol[1]{%
   \ifcase#1\or *\or **\or \dagger\or \ddagger\or
   \mathsection\or \mathparagraph\or \|\or \dagger\dagger \or \ddagger\ddagger \fi}
\makeatother

\begin{document}
\begin{frontmatter}

\title{Accelerated ADMM: Automated Parameter Tuning and Improved Linear Convergence}

\author[First,Third]{M. Tavakoli}
\author[Second,Third]{F. Jakob}
\author[First]{G. Carnevale}
\author[First]{G. Notarstefano}
\author[Second]{A. Iannelli}

\thanks[footnoteinfo]{F. Jakob acknowledges the support of the International Max Planck Research School for Intelligent Systems (IMPRS-IS).}

\address[First]{Università di Bologna, Department of Electrical, Electronic and Information
Engineering, Bologna, Italy 
(e-mail: meisam.tavakoli@studio.unibo.it, guido.carnevale, giuseppe.notarstefano@unibo.it).}

\address[Second]{University of Stuttgart, Institute for Systems Theory and Automatic Control, Stuttgart, Germany 
(e-mail: fabian.jakob, andrea.iannelli@ist.uni-stuttgart.de).}

\address[Third]{The first two authors contributed equally to this work.}

\begin{abstract}                % Abstract of 50--100 words
This work studies the linear convergence of an accelerated scheme of the Alternating Direction Method of Multipliers (ADMM) for strongly convex and Lipschitz-smooth problems. We use the methodology of expressing the accelerated ADMM as a Lur'e system, i.e., an interconnection of a linear dynamical system in feedback with a slope-restricted operator, and we use Integral Quadratic Constraints to establish linear convergence. We leverage this machinery to systematically explore parameter tuning heuristics, including Nesterov-inspired choices and configurations identified via grid search, and analyze their impact on the convergence rate. Our new bounds show improved linear convergence rates compared to the vanilla algorithm and previously proposed accelerated variants, which is also empirically validated on a LASSO regression benchmark.
\end{abstract}


\end{frontmatter}
%===============================================================================

\section{Introduction}
\vspace{-3mm}
The Alternating Direction Method of Multipliers (ADMM) \citep{Glowinski1975, boydAdmm} has emerged as a benchmark
algorithm 
for solving convex composite optimization problems.
Its ability to decompose the original problem into subproblems that can be solved in parallel made it especially popular in computation-heavy domains such as signal processing, machine learning, or distributed control.

\meisam{It has long been established} that ADMM enjoys provable convergence guarantees in the convex 
setting. \cite{DengYin} showed that for general convex problems, ADMM 
achieves an $\mathcal{O}(1/k)$ convergence rate, while for the
strongly convex and smooth case it attains a linear convergence rate, with optimally tuned parameters given in \cite{pontus_boyd}.
An alternative proof of linear convergence was presented by~\cite{NishiharaADMM}, who applied the Integral Quadratic Constraints (IQC) framework \citep{Megretski1997,LessardRechtPackard,SchererEbenbauer} to numerically derive convergence rates and tune parameters. The proof frames the optimization algorithm as a Lur'e system (a feedback interconnection of an LTI system and monotone operators), and casts the convergence analysis as a stability problem. This framework offers the main advantage that the analysis of altered schemes, such as over-relaxed variants, can also be seamlessly carried out in an automated fashion as long as they can be framed in the Lur'e setting.

Accelerated variants of ADMM (A-ADMM) have also been proposed and analyzed in recent years.
\cite{FastADMM} showed that Nesterov acceleration leads to a
convergence rate of $\mathcal{O}(1/k^2)$ in the strongly convex setting. \cite{Sterck} provided linear convergence guarantees when Anderson acceleration is applied.
\cite{PatrinosStellaBemporad, PejcicJones} proved linear convergence of accelerated Douglas–Rachford splitting and ADMM for strongly convex quadratic problems, showing that the optimal rate matches that of Nesterov’s fast gradient method. However, an analysis of the convergence rates of A-ADMM for general strongly convex and smooth problems is still unaddressed, as well as the question of whether other known rates from classical accelerated gradient methods (AGMs) can be matched.

\begin{table*}
\centering
\caption{Contextualization of the results. Here,
$\rho_1 = 1 - \frac{1}{\sqrt{\kappa}}$,
$\rho_2 = \sqrt{1-\frac{\sqrt{2\kappa - 1}}{\kappa}}$, 
with $\kappa = \frac{L}{m}$. NM stands for Nesterov parameter selection (Table~\ref{tb:param_selection}).}
\begin{tabular}{ccccc}
\toprule
Setting & Objective $f$ & Algorithm & Rate $\rho$ & Reference \\
\midrule
$\min f(x)$ & Quadr. & Nesterov Method & $\rho_1$ & \cite{Nesterov2004} \\
$\min f(x)$ & $S(m,L)$ & Nesterov Method & $\rho_2$ & \cite{safavi} \\
\eqref{constrained optimization} & Quadr. & A-ADMM (NM) & $\rho_1$ & \cite{PatrinosStellaBemporad,PejcicJones} \\
\eqref{constrained optimization} & $S(m,L)$ & A-ADMM (NM) & Fig.~\ref{fig:NM_study} & This work \\
\bottomrule
\end{tabular}
\label{tab:context}
\end{table*}

In this work, we extend the IQC-based analysis of vanilla ADMM \citep{NishiharaADMM} to the accelerated case, when one of the objectives is strongly convex and smooth. By casting the A-ADMM algorithm as a Lur'e system, we propose a semi-definite program (SDP) that provides a numerical tool to theoretically verify worst-case convergence rates of accelerated schemes \meisam{using different parameter selections}. 
\meisam{Crucially, we show that dynamic O'Shea-Zames-Falb (OZF) IQCs \citep{Scherer} are essential to certify convergence of the accelerated scheme. Based on this, we propose a systematic parameter tuning procedure and demonstrate that Nesterov-inspired parameter selections match the rates of the analogous AGM in the strongly convex and smooth setting
(cf. Table~\ref{tab:context}). }
Finally, we identify through grid search a new A-ADMM configuration with the fastest certified rate. Our theoretical results are validated on a LASSO regression case study, demonstrating consistently improved convergence speed of our schemes compared to existing benchmark algorithms. 

The remainder of the paper is organized as follows. Section~2 presents the problem setup, 
the ADMM formulation, and the required IQC preliminaries. Section~3 provides a Lur'e representation of the A-ADMM scheme. In Section~4, we study the convergence rates of this A-ADMM scheme under different parameter configurations. In Section~5, we
provide a case study on LASSO regression. Finally,
Section~6 concludes the paper.

\textbf{Notation}. The identity matrix of dimension $p$ is denoted as $I_p$. The gradient and subdifferential of a function $f$ are denoted by
$\nabla f$ and $\partial f$, respectively. For $0 < m \leq L < \infty$, we let $S_p(m,L)$
denote the class of functions $f:\mathbb{R}^p \to \mathbb{R}$ that
are $m$-strongly convex and have $L$-Lipschitz continuous gradients. The special case
$S_p(0,\infty)$ corresponds to the set of proper, closed and convex functions. For a matrix $M$,
its condition number is defined as $\kappa_M = \bar{\sigma}(M)/\underline{\sigma}(M)$,
where $\bar{\sigma}(M)$ and $\underline{\sigma}(M)$ denote the largest and smallest
singular values of $M$, respectively. A discrete-time LTI system maps the signals $u\mapsto y$ via the recursion $\xi_{k+1} = A\xi_k + Bu_k, \; y_k = C\xi_k + Du_k$ as a $\xi_0$ dependent mapping, and will be compactly expressed as  
$y = G u$ with $G={\small \left[\begin{array}{c|c}
A & B \\[-1pt]\hline
C & D 
\end{array}\right]}$. We will write $G \otimes I_d$ to define a lifted system {\small $\left[\begin{array}{c|c}
A \!\otimes\! I_d & B \!\otimes\! I_d \\[-1pt]\hline
C \!\otimes\! I_d & D \!\otimes\! I_d
\end{array}\right]$}. The series interconnection between two LTI systems $G_1: u\mapsto y$ and $G_2: y \mapsto w$ will be denoted as $G_1 \cdot G_2: u \mapsto w$. The forward shift operator is denoted as $\mathbf{z}$. 

\vspace{-1mm}
\section{PRELIMINARIES}
\label{sec:PRELIMINARIES}
\vspace{-2mm}
\subsection{Accelerated ADMM}
\vspace{-3mm}
We consider constrained convex optimization problems of the form
\begin{equation}
    \min_{x \in \mathbb{R}^p, \; z \in \mathbb{R}^q} \; f(x) + g(z) \quad \text{s.t.} \quad Ax + Bz = c,
    \label{constrained optimization}
\end{equation}
where \( f: \mathbb{R}^p \to \mathbb{R} \), $f \in S_p(m,L)$ and \( g: \mathbb{R}^q \to \mathbb{R} \), $g \in S_q(0,\infty)$. The matrices \( A \in \mathbb{R}^{p
\times p} \), \( B \in \mathbb{R}^{p \times q} \), and \( c \in \mathbb{R}^p \) define
linear equality constraints coupling the variables \( x \) and \( z \), where we
particularly assume that $A$ is invertible and that $B$ has full column rank. Having a
strongly convex and smooth component in~\eqref{constrained optimization} is common in many practical
applications~\citep{DengYin}, while restricting $A$ to be square and invertible is slightly restrictive, but also done in \cite{NishiharaADMM}. We note that many practical problem instances \meisam{still} satisfy this structure, including consensus problems and various distributed optimization tasks \citep{NotarstefanoNotarnicolaCamisa}.

To address problem~\eqref{constrained optimization}, we use the over-relaxed ADMM algorithm (OR-ADMM) \citep{boydAdmm} 
\begin{subequations}\label{eq:ADMM_whole}
\begin{align}
    \label{eq:ADMM_whole:x}
    x_{k+1} &= \arg\min_{x\in\mathbb{R}^p} \; f(x) + \tfrac{1}{2\nu_1}\,\|Ax + B\hat{z}_k - c + \hat{\lambda}_k\|^2 \\
    \label{eq:ADMM_whole:z}
    z_{k+1} &= \arg\min_{z\in\mathbb{R}^q} \; g(z) + \tfrac{1}{2\nu_1}\,\|\alpha Ax_{k+1} - (1 - \alpha) B\hat z_k  \\[-1ex]
&\qquad \qquad \qquad \qquad \qquad \qquad \qquad + Bz - \alpha c + \hat{\lambda}_k\|^2 \notag \\
    \lambda_{k+1} &= \hat{\lambda}_k + \alpha Ax_{k+1} - (1 - \alpha) B\hat z_k + Bz_{k+1} - \alpha c
\end{align}
\end{subequations}
as a starting point, and augment it with a \emph{momentum term}
\vspace{-3mm}
\begin{subequations}\label{eq:ADMM_accel}
\begin{align}  
    \label{eq:ADMM_accel:z} 
    \hat{z}_k &= z_k + \nu_2 (z_k - z_{k-1}) \\
    \label{eq:ADMM_accel:lambda} 
    \hat{\lambda}_k &= \lambda_k + \nu_2 (\lambda_k - \lambda_{k-1}),
\end{align}
\end{subequations}
as in \cite{FastADMM,PejcicJones}.
Here, $\nu_1 > 0$ is a % step size
\meisam{penalty parameter}, $\nu_2 \geq 0$ a momentum parameter, and $\alpha\geq 0$ the so-called over-relaxation parameter. 
We recover the nominal A-ADMM by setting $\alpha=1$, and the vanilla, non-accelerated, ADMM with $\nu_2 = 0$.
Both $\alpha$ and $\nu_2$ have been shown to improve the convergence speed for suitable choices \citep{Eckstein1992,PatrinosStellaBemporad}. The over-relaxation parameter $\alpha$ is chosen in the interval $(0,2]$, while $\nu_2$ is often selected according to Nesterov's schemes \citep{Nesterov2004}, e.g., as
in \cite{PatrinosStellaBemporad}. 
To explore alternative parameter selection strategies, we employ the IQC framework.
\vspace{-2mm}
%%%%%%%%%%%%%%%%%%%%%%%%%%%%%%%%%%%%%%%%%%%%%
\subsection{Integral Quadratic Constraints for Convex Functions} \label{IQC-section}
\vspace{-3mm}
Algorithm \eqref{eq:ADMM_whole}-\eqref{eq:ADMM_accel} will be analyzed as a dynamical system converging to points that satisfy the first-order optimality conditions of \eqref{constrained optimization}. We will make use of the well-known fact that
(sub)gradients of (strongly) convex functions are monotone and, thus, their input/output relation satisfies IQCs for slope-restricted operators \citep{LessardRechtPackard}. 
To be more precise, they have been shown to fulfill so-called $\rho$-hard O'Shea-Zames–Falb IQCs.

\begin{prop}[\cite{Scherer}]\label{prop:ZF_IQC}
    Let $f \in S_p(m,L)$, $g \in S_p(0,\infty)$ and let $\rho\in(0,1)$ be an exponential discount factor. For $\bar{n} \in \mathbb{N}_0$, let $\{h^f_\tau\}_{\rev{\tau=- \bar n}}^{\bar{n}}$ and $\{h^g_\tau\}_{\rev{\tau=- \bar n}}^{\bar{n}}$ be sequences of filter coefficients satisfying \rev{$h_0^f, h_0^g \geq 0$ and
    \begin{equation}\label{eq:filter_conditions}
        h_\tau^* \leq 0 \,\,\, \forall \tau\neq 0, \,\,\, \sum_{\tau=-\bar{n}}^{\bar{n}} \rho^{-2\tau} h_\tau^* \geq 0, \,\,\, \sum_{\tau=-\bar{n}}^{\bar{n}} \rho^{2\tau} h_\tau^* \geq 0
    \end{equation}
    for $* \in \{f,g\}$},
    and define the parametrized LTI filters
    \begin{align}
    \begin{aligned}
    \Psi_g(h^g)&=
        {
            \begin{bmatrix}
            {\tiny \displaystyle \sum_{\tau = 0}^{\bar{n}}} h_\tau^g \mathbf{z} ^{-\tau} & 0 \\[-8pt] 0 & \rev{{\tiny \displaystyle \sum_{\tau = -\bar{n}}^{0}} h_{\tau}^g \mathbf{z} ^{\tau}}
            \end{bmatrix}
        }
        \\ 
        \Psi_f(h^f) &= \Psi_g(h^f) 
        % {\setlength{\arraycolsep}{1pt}
            \begin{bmatrix}L&-1\\-m&1\end{bmatrix}
        % }.
    \end{aligned}
    \end{align}
    Let $\{a_k\}, \{b_k\}$ be some $p$-dim. square summable sequences and let $a^\star, b^\star \in \mathbb{R}^p$ be constant references. For $\gamma_k \in \partial g(b_k)$ and $\gamma^\star \in \partial g(b^\star)$, define
    \begin{align*}
        \tilde{a}_k &= a_k - a^\star, &  {\nabla \tilde f}_k &= \nabla f(a_k) - \nabla f(a^\star), \\
        \tilde{b}_k &= b_k - b^\star, & \tilde{\gamma}_k &= \gamma_k - \gamma^\star.
    \end{align*}
    Finally, define the filtered sequences
    \begin{equation*}
        \psi_{1,k} = (\Psi_f \otimes I_p) \begin{bmatrix}
            \tilde{a}_k \\ {\nabla \tilde f}_k 
        \end{bmatrix}, \quad
        \psi_{2,k} = (\Psi_g \otimes I_p) \begin{bmatrix}
            \tilde{b}_k \\ \tilde \gamma_k
        \end{bmatrix}
    \end{equation*}
    and $M = \left[\begin{smallmatrix}
        0 & 1 \\ 1  & 0
    \end{smallmatrix}\right]$. Then, for all $T\geq 0$, it holds
    \begin{align}\label{eq:IQC_quad}
        \begin{aligned}
        &\sum_{k=0}^T \rho^{-2k} \psi_{i,k}^\top (M\otimes I_p) \psi_{i,k} \geq 0, \quad i=1,2.
        \end{aligned}
    \end{align}
\end{prop}
\vspace{-2mm}
The filter coefficients $h^f$, $h^g$ 
will serve as degrees-of-freedom in the resulting SDP, subject to the convex constraint \eqref{eq:filter_conditions}. In theory, the filter dimension $\bar{n}$ can be infinite, but is chosen finite in practice to yield finite-dimensional state-space realizations. A choice of $\bar{n} = 0$ leads to static filters $\Psi_f,\Psi_g$ and pointwise satisfaction of the inequalities \eqref{eq:IQC_quad}, which is often referred to in the literature as a sector IQC \citep{LessardRechtPackard}. \meisam{While sector IQCs are sufficient to certify the convergence of vanilla ADMM
\citep{NishiharaADMM}, we will show in Section~\ref{sec:param_selection} that the more general
dynamic OZF IQCs $(\bar{n} > 0)$ are essential for the accelerated case.}

\vspace{-2mm}
%%%%%%%%%%%%%%%%%%%%%%%%%%%%%%%%%%%%%%%%%%%%%%%%%%%%%%%%%%%%%%%%%
\subsection{Stability Analysis of First-Order Algorithms}
\vspace{-3mm}
Consider the following Lur'e representation of a first-order algorithm
\begin{subequations}\label{eq:Lure}
\begin{align}\label{eq:LTI_sys}
    \begin{aligned}
    \xi_{k+1} &= (\hat{A} \otimes I_p) \xi_k + (\hat{B} \otimes I_p) w_k \\
    \hat{v}_k &= (\hat{C} \otimes I_p) \xi_k + (\hat{D} \otimes I_p) w_k
    \end{aligned}
\end{align}
with state $\xi_k \in \mathbb{R}^{n_\xi}$, output $\hat{v}_k = v_k+\bar{v}$ for some constant offset $\bar{v}$, and
\begin{equation}
 v_k = \begin{bmatrix} v_{1,k} \\ v_{2,k} \end{bmatrix}, \quad
    w_k = \begin{bmatrix} \nabla \hat{f}(v_{1,k}) \\ \gamma_k \end{bmatrix}, \, \gamma_k \in \partial \hat{g}(v_{2,k}),
\end{equation}
\end{subequations}
where $\hat{f}\in S_p(\hat{m},\hat{L})$, $\hat{g}\in S_p(0,\infty)$ for some $\hat{m},\hat{L}$.
It has been shown that ADMM, as well as many other composite optimization algorithms, can be represented as \eqref{eq:Lure} 
\citep{NishiharaADMM} when $\hat{f}, \hat{g}$ are suitably chosen.
In general, $(\hat A,\hat B,\hat C,\hat D)$ fulfills structural assumptions such that its fixed-point $(\xi^\star, \hat{v}^\star, w^\star)$ is unique and satisfies first-order optimality \citep{Upadhyaya2024}. Showing exponential stability of \eqref{eq:Lure} is therefore the same as showing linear convergence of the underlying algorithm.
With a coordinate shift $(\tilde \xi_k, \tilde{v}_k, \tilde w_k) \triangleq (\xi_k-\xi^\star, \hat{v}_k-\hat{v}^\star, w_k-w^\star)$ it is straightforward to show that the error coordinates evolve with the same state-space description \eqref{eq:LTI_sys}, i.e.,
\begin{equation}
    \tilde{v} = \Biggl( \underbrace{\left[\begin{array}{c|c}  \hat A &  \hat B \\ \hline  \hat C &  \hat D \end{array}\right]}_{=: G} \otimes I_p \Biggr) \tilde{w}.
\end{equation}
Note that $\tilde{v}_k = \hat v_k - \hat v^\star = v_k - v^\star$ (eliminating the constant offset), so that by Proposition~\ref{prop:ZF_IQC} the gradient and subgradient components of $(\tilde{v},\tilde{w})$ satisfy a $\rho$-OZF IQC. With a suitable permutation and stacking of the filters $\Psi_1, \Psi_2$, we can form a compact filter $\Psi$ such that 
\begin{equation}
\begin{aligned}
    \psi_k \triangleq 
    \begin{bmatrix} \psi_{1,k} \\ \psi_{2,k} \end{bmatrix}
    = (\Psi \otimes I_p) \begin{bmatrix} \tilde v_k \\ \tilde  w_k \end{bmatrix},
\end{aligned}
\end{equation}
with $\psi_{1,k}, \psi_{2,k}$ as in Proposition~\ref{prop:ZF_IQC} (with $a\triangleq v_1,\, b\triangleq v_2$).
Define the augmented plant as the mapping $\tilde w \mapsto \psi$, realized as
\begin{equation}\label{eq:augmented_plant}
\begin{aligned}
    \Psi \cdot \begin{bmatrix} G \\ I_2 \end{bmatrix} 
    \triangleq \left[\begin{array}{c|c} \mathbf{A} & \mathbf B \\ \hline \mathbf C & \mathbf D \end{array}\right],
\end{aligned}
\end{equation}
\rev{
where $\left[ \begin{smallmatrix} G \\ I_2 \end{smallmatrix} \right]$ denotes the dynamical system arising by vertically concatenating the transfer matrix of $G$ and a static feedthrough, mapping $\tilde{w} \to \left[ \begin{smallmatrix} \tilde{v} \\ \tilde{w} \end{smallmatrix} \right]$.}
Then the following convergence theorem holds.

\begin{thm}\label{thm:worst-case cng thm}
Take $\rho \in(0,1)$, and filters $\Psi_1, \Psi_2$ that satisfy Proposition~\ref{prop:ZF_IQC} with $m= \hat{m}, L=\hat{L}$. Form the augmented plant \eqref{eq:augmented_plant}.
If there exist $P\succ0$ and filter coefficients of $\Psi_1,\Psi_2$ such that the matrix inequality
\begin{equation}\label{eq:lmi_condition}
\begin{bmatrix}
\mathbf A^\top P \mathbf A - \rho^2 P & \mathbf A^\top P \mathbf B \\
\mathbf B^\top P \mathbf A            & \mathbf B^\top P \mathbf B
\end{bmatrix}
+
\begin{bmatrix}\mathbf C^\top\\ \mathbf D^\top\end{bmatrix}
\begin{bmatrix}
    M & 0 \\ 0 & M
\end{bmatrix}
\begin{bmatrix}\mathbf C & \mathbf D\end{bmatrix}
\;\preceq\;0
\end{equation}
holds, then \eqref{eq:Lure} is exponentially stable with rate $\rho$, i.e.,
\begin{equation}
    \| \xi_k - \xi^\star \| \leq \sqrt{\kappa_P}\, \rho^k \| \xi_0 - \xi^\star \|.
\end{equation}
\end{thm}
The proof comes as a straightforward extension of \citep[Theorem~4]{LessardRechtPackard} to the two-operator
case, which particularly exploits the positivity constraint \eqref{eq:IQC_quad}.
Eq. \eqref{eq:lmi_condition} is independent of the problem dimension $p$, and can be used to determine the minimum worst-case convergence rate $\rho$ via bisection.
Note that Theorem~\ref{thm:worst-case cng thm} recovers \cite[Theorem~6]{NishiharaADMM}
for $\bar{n}=0$.
\vspace{-3mm}
\section{A-ADMM as a Dynamical System} \label{sec:main results}
\vspace{-3mm}
Building on the previous derivations, we now formulate the A-ADMM \eqref{eq:ADMM_whole}-\eqref{eq:ADMM_accel} as a Lur'e system \eqref{eq:Lure}. In line with \cite{NishiharaADMM}, we introduce the coordinate change \( r_k
= Ax_k \), \( s_k = Bz_k\). Moreover, define
\begin{equation}
    \hat{f} = f \circ A^{-1}, \quad \hat{g} = g \circ B^\dagger + \mathcal{I}_{\mathrm{im}\,B},
\end{equation}
where $B^\dagger$ is
a left-inverse of $B$ and $\mathcal{I}_{\mathrm{im}\,B}$ is the indicator of the image of
$B$. It is straightforward to verify that $\hat{g}\in S_p(0, \infty)$ and $\hat{f}\in S_p(\hat m, \hat L)$
with
\begin{align} \label{eq:normalizing}
    \hat{m} = \tfrac{m}{\bar{\sigma}^2(A)}, \quad 
    \hat{L} = \tfrac{L}{\underline{\sigma}^2(A)}.
\end{align}
Accordingly, we define the condition number of the problem as $\kappa = \frac{\hat{L}}{\hat m} =\kappa_f \kappa_A^2$.

Consequently, the updates \eqref{eq:ADMM_whole:x}, \eqref{eq:ADMM_whole:z}, \eqref{eq:ADMM_accel:z} can be rewritten as
\begin{subequations}
\begin{align}
    x_{k+1} &= A^{-1} \arg\min_{r\in\mathbb{R}^p}\; \hat f(r) 
    + \tfrac{1}{2\nu_1}\,\lVert r + \hat s_k - c + \hat \lambda_k \rVert^2 \\[-1ex]
    z_{k+1} &= B^\dagger \arg\min_{s\in\mathbb{R}^p}\;  \hat g(s)
    + \tfrac{1}{2\nu_1}\,\lVert \alpha r_{k+1} - \\[-1ex]
    \notag
     &\hspace{35mm} (1-\alpha) \hat s_k + s - \alpha c + \hat \lambda_k \rVert^2 \\[0.2em]
     \label{eq:ADMM_transformed:accel:s}
     \hat{s}_k &= s_k + \nu_2(s_k - s_{k-1}).
\end{align}
\end{subequations}
Using the definition of the proximal operator
\begin{align*}
\mathrm{prox}_f(z) := \arg\min_{x} f(x) + \frac{1}{2} \|x - z\|^2 ,
\end{align*}
we can summarize the A-ADMM in transformed coordinates compactly as 
\begin{subequations}\label{eq:ADMM_transformed}
\begin{align}
    r_{k+1} &= \operatorname{prox}_{\nu_1 \hat f}\!\big(c - \hat{s}_k - \hat{\lambda}_k\big) \label{eq:ADMM:r-update} \\
    s_{k+1} &= \operatorname{prox}_{\nu_1 \hat g}\!\big(\alpha c - \alpha r_{k+1} - (\alpha-1)\hat{s}_k - \hat{\lambda}_k\big) \label{eq:ADMM:s-update} \\ 
     \lambda_{k+1} &= \alpha r_{k+1} + (\alpha-1)\hat{s}_k + s_{k+1} - \alpha c + \hat{\lambda}_k,
\end{align}
\end{subequations}
with $\hat{s}_k, \hat{\lambda}_k$ as in \eqref{eq:ADMM_transformed:accel:s}, \eqref{eq:ADMM_accel:lambda}.
Next, we bring the recursion \eqref{eq:ADMM_transformed} into the form \eqref{eq:Lure}.

\begin{prop}\label{prop:a-admm}
Define the state $\xi_k \triangleq \begin{bmatrix}
    \lambda_{k-1}^\top & s_{k-1}^\top & \lambda_{k}^\top & s_{k}^\top
\end{bmatrix}^\top$, output $v_k \triangleq \begin{bmatrix}
    r_{k+1}^\top & s_{k+1}^\top
\end{bmatrix}^\top$, offset $\bar{v}\triangleq \begin{bmatrix}
    -c^\top & 0
\end{bmatrix}^\top$, and input
$w_k \triangleq \begin{bmatrix}
    \nabla \hat{f}(r_{k+1})^\top & \gamma_k^\top
\end{bmatrix}^\top$ for any 
$\gamma_k \in \partial \hat{g}(s_{k+1})$. Then the sequences $\xi_k$, $w_k$, and
$\hat v_k=v_k + \bar{v}$ \rev{generated by the A\nobreakdash-ADMM recursion \eqref{eq:ADMM_transformed}} satisfy~\eqref{eq:Lure} with the matrices
\begin{align}\label{prop-OR-A-ADMM-ss1}
\begin{aligned}
\hat A &= 
\begin{pmatrix}
0 & 0 & 1 & 0 \\
0 & 0 & 0 & 1 \\
0 & 0 & 0 & 0 \\
-\nu_2(\alpha - 1) & -\nu_2 & (\alpha - 1)(1+\nu_2) & 1+\nu_2
\end{pmatrix}, \\[1ex]
\hat C &=
\begin{pmatrix}
\nu_2 & \nu_2 & -(1+\nu_2) & -(1+\nu_2)\\
-\nu_2(\alpha - 1) & -\nu_2 & (\alpha - 1)(1+\nu_2) & 1+\nu_2
\end{pmatrix}, \\
\hat B &=
\begin{pmatrix}
0 & 0 \\
0 & 0 \\
0 & -\nu_1 \\
\alpha \nu_1 & -\nu_1
\end{pmatrix},\quad
\hat D =
\begin{pmatrix}
 -\nu_1 & 0\\
 \alpha \nu_1 & -\nu_1
\end{pmatrix}.
\end{aligned}
\end{align}
\label{prop:state and input}
\end{prop}
\vspace{-5mm}
\begin{pf}
Apply the first-order optimality condition to the $r_{k+1}$ and $s_{k+1}$ updates \eqref{eq:ADMM:r-update}-\eqref{eq:ADMM:s-update}
\begin{subequations}\label{eq:first order nc}
\begin{align}
    0 &= \nu_1 \nabla \hat{f}(r_{k+1}) + r_{k+1} - c + \hat{s}_k + \hat{\lambda}_k,\\[0.5ex]
    0 &\in \nu_1 \partial \hat{g}(s_{k+1}) + s_{k+1} - \alpha c + \alpha r_{k+1} + (\alpha-1)\hat s_k + \hat \lambda_k.
\end{align}
\end{subequations} 
Take some $\gamma_k \in \partial \hat g(s_{k+1})$, solve for $r_{k+1}$ and $s_{k+1}$, and plug \meisam{them} into the dual update $\lambda_{k+1}$
to get
\begin{subequations} \label{eq:first order nc2}
\begin{align}  
    r_{k+1} &= -\hat{s}_k - \hat{\lambda}_k + c -\nu_1 \nabla \hat{f}(r_{k+1}),
    \label{eq:fnc_r} \\[0.5ex]
    s_{k+1} &= \hat s_k + (\alpha-1) \hat{\lambda}_k  + \alpha \nu_1 \nabla \hat f(r_{k+1}) - \nu_1 \gamma_k, \label{eq:fnc_g} \\
    \lambda_{k+1} &=  - \nu_1 \gamma_k. \label{eq:fnc_l}
\end{align}
\end{subequations}
Now plug in the definitions of $\hat s_k, \hat \lambda_k$ and observe that $s_{k+1}, \lambda_{k+1}$, and $r_{k+1}-c$ can all be written as linear combinations of the state $\xi_k$ and input $w_k$. It is then straightforward to bring $\xi_{k+1}$ and $\hat v_k$ into the matrix form \eqref{eq:Lure}.
\hfill\qed
\end{pf}

\rev{Theorem \ref{thm:worst-case cng thm} together with Proposition \ref{prop:a-admm} provides a numerical tool to certify upper bounds on the worst-case convergence rate of A-ADMM. Next, we explore the influence of acceleration and the algorithm parameters on these rates.}
\vspace{-3mm}
%%%%%%%%%%%%%%%%%%%%%%%%%%%%%%%%%%%%%%%%%%%%%%%%%%%%%%%%%%%%%%%%%%%%%%%
\section{Parameter Selection}
\label{sec:param_selection} 
\vspace{-3mm}

\rev{We investigate different parameter selection heuristics for $\nu_1$ and $\nu_2$.}
As initial heuristics, we adopt the same parameter tuning schemes of \meisam{AGMs.} Next, we perform a grid search to identify potentially superior configurations. 
\rev{
\meisam{As a baseline for comparison}, we use a vanilla ADMM instance with the parameter $\nu_1$ tuned according to the optimal choice in the strongly convex and smooth setting characterized by~\cite{pontus_boyd}.}
\vspace{-2mm}
\subsection{AGM-inspired Parameters}
\vspace{-2mm}
\rev{
Since Nesterov-inspired (NM) tuning (cf. Table~\ref{tb:param_selection}) \meisam{is a standard benchmark among AGMs} and has proven effective for quadratic ADMM problems \citep{PatrinosStellaBemporad}, it is also a natural candidate for the strongly-convex-smooth ADMM.}
We use the machinery developed in the previous sections to numerically verify the convergence rates $\rho$ for this choice of $(\nu_1,\nu_2)$.
\rev{By bisection, we find the smallest $\rho$ rendering \eqref{eq:lmi_condition} feasible} and plot its value as a function of the condition number $\kappa = \hat L / \hat m$\footnote{Code for all experiments can be found on \\
https://github.com/col-tasas/2025-accelerated-admm
}.

\rev{
The resulting convergence rate curve is reported in Fig.~\ref{fig:NM_study}.
\meisam{To highlight the role of dynamic IQC multipliers, we plot the curve resulting from different choices of the OZF order~$\bar{n}$}.
We emphasize that the choice of a sector IQC ($\bar{n}=0$) is not sufficient to certify convergence for all condition numbers $\kappa$. This observation is consistent with AGMs \citep{LessardRechtPackard}. Increasing the OZF order to $\bar{n}=4$ improves the rate substantially, which then closely tracks the classical NM rate $\rho = {\scriptsize\sqrt{1-\tfrac{\sqrt{2\kappa - 1}}{\kappa}}}$ \citep{safavi}. The rate is matched for condition ratios larger than $\kappa \approx 10$; with deviations for low $\kappa$ that are also typical of vanilla ADMM \citep{NishiharaADMM}.
}

\rev{
Our framework also allows us to seamlessly evaluate alternative parameters. For instance, since the Triple Momentum (TM) algorithm \citep{Freeman2028} is known to outperform NM in the standard strongly-convex-smooth setting, we also investigate its application to A-ADMM. The comparison curves for vanilla ADMM, and A-ADMM with NM- and TM-tuning (cf. Table~\ref{tb:param_selection}) are reported in Fig.~\ref{fig:A-ADMM}. We observe that while TM parameters offer advantages in certain $\kappa$ regimes, the certified rate exceeds one for high condition numbers, even with high-order OZF multipliers. Thus, this configuration does not yet provide a globally certified improvement over the NM scheme in the ADMM setting.
}

\begin{table}
\centering
\caption{\rev{AGM-inspired} tuning of A-ADMM, where $\rho:=1 - 1/\sqrt{\kappa}$.}
\label{tb:param_selection}
\begin{tabular}{lcc}
\toprule
Method & $\nu_1$ & $\nu_2$ \\ \midrule
NM  & $\frac{1}{\hat{L}}$ & $\frac{\sqrt{\hat{L}} - \sqrt{\hat{m}}}{\sqrt{\hat{L}} + \sqrt{\hat{m}}}$ \\
TM & $\frac{1+\rho}{\hat L}$ & $\frac{\rho^2}{2-\rho}$\\ \bottomrule
\end{tabular}
\end{table}

\begin{figure}
    \centering
        \centering
        % This file was created with tikzplotlib v0.10.1.
\begin{tikzpicture}

\definecolor{darkgray176}{RGB}{176,176,176}
\definecolor{mycolor1}{RGB}{204,51,63}   % red
\definecolor{mycolor2}{RGB}{68,12,84}    % purple
\definecolor{mycolor3}{RGB}{53,95,141}   % blue
\definecolor{mycolor4}{RGB}{34,168,132}  % teal

\definecolor{ozf0}{RGB}{130,202,225}
\definecolor{ozf1}{RGB}{66,146,198}
\definecolor{ozf4}{RGB}{8,48,107}
\definecolor{refgray}{RGB}{110,110,110}

\begin{axis}[
    width=7cm,
    height=4cm,
    % every axis/.append style={font=\small},
    % tick label style={font=\small},
    legend style={
      fill opacity=0.4,
      draw opacity=1,
      text opacity=1,
      font=\small,
      nodes={scale=0.7, transform shape},
      at={(0.97,0.03)},
      anchor=south east,
      legend cell align={left}
    },
    log basis x={10},
    tick align=inside,
    tick pos=left,
    x grid style={darkgray176},
    xlabel={$\kappa$},
    xlabel style={
      at={(axis description cs:0.5,-0.04)},
      anchor=north
    },
    xmajorgrids,
    xmin=0.707945784384138, xmax=1412.53754462275,
    xmode=log,
    xtick style={color=black},
    xtick={0.01,0.1,1,10,100,1000,10000,100000,1000000},
    xticklabels={
      \(\displaystyle {10^{-2}}\),
      \(\displaystyle {10^{-1}}\),
      \(\displaystyle {10^{0}}\),
      \(\displaystyle {10^{1}}\),
      \(\displaystyle {10^{2}}\),
      \(\displaystyle {10^{3}}\),
      \(\displaystyle {10^{4}}\),
      \(\displaystyle {10^{5}}\),
      \(\displaystyle {10^{6}}\)
    },
    y grid style={darkgray176},
    ylabel={$\rho$},
    ymajorgrids,
    ymin=0.45, ymax=1.08875,
    ytick style={color=black, yshift=-3mm}
]

% --- ADMM ---
% \addplot [mycolor1, thick, mark=*, dashed, mark size=1.25pt, mark options={solid, fill=mycolor1, draw=mycolor1}]
% table {%
% 1.00 0.5002
% 1.44 0.5440
% 2.07 0.5891
% 2.98 0.6322
% 4.28 0.6735
% 6.16 0.7122
% 8.86 0.7484
% 12.74 0.7814
% 18.33 0.8106
% 26.37 0.8373
% 37.93 0.8607
% 54.56 0.8811
% 78.48 0.8988
% 112.88 0.9141
% 162.38 0.9274
% 233.57 0.9388
% 335.98 0.9483
% 483.29 0.9566
% 695.19 0.9636
% 1000.00 0.9699
% };
% \addlegendentry{ADMM}

% --- A-ADMM (n_OZF=0) ---
% \addplot [mycolor1, thick, mark=*, dashed, mark size=1.5pt, mark options={solid, fill=mycolor1, draw=mycolor1}]
% table {%
% 1 0.5001953125
% 1.27427498570313 0.5560546875
% 1.62377673918872 0.61318359375
% 2.06913808111479 0.669677734375
% 2.63665089873036 0.7236328125
% 3.35981828628378 0.77568359375
% 4.28133239871939 0.835986328125
% 5.45559478116852 0.896923828125
% 6.95192796177561 0.952783203125
% 8.85866790410082 1.003564453125
% 11.2883789168469 1.0486328125 
% 14.3844988828766 1.088623046875 
% 18.3298071083244 1.124169921875 
% 23.3572146909012 1.1552734375 
% 29.7635144163132 1.183203125 
% 37.9269019073225 1.207958984375 
% 48.3293023857175 1.23017578125 
% 61.5848211066026 1.249853515625 
% 78.4759970351461 1.267626953125 
% 100 1.28349609375
% };
% \addlegendentry{A-ADMM ($n_{\mathrm{OZF}}=0$)}

% --- A-ADMM (n_OZF=0) ---
\addplot [ozf0, thick, mark=*, dashed, mark size=1.25pt, mark options={solid, fill=ozf0, draw=ozf0}]
table {%
1 0.5002
1.44 0.5434
2.07 0.5872
2.98 0.6329
4.28 0.6735
6.16 0.7217
8.86 0.7871
12.74 0.8480
18.33 0.9014
26.37 0.9464
37.93 0.9852
54.56 1.0169
78.48 1.0442
112.88 1.0664
% 162.38 1.0848
% 233.57 1.1058
% 335.98 1.1299
% 483.29 1.1527
% 695.19 1.1699
% 1000 1.1838
};
\addlegendentry{A-ADMM ($\bar n = {0}$)}

% --- A-ADMM (n_OZF=1) ---
\addplot [ozf1, thick, mark=*, dashed, mark size=1.25pt, mark options={solid, fill=ozf1, draw=ozf1}]
table {%
1 0.5002
1.44 0.5415
2.07 0.5878
2.98 0.6316
4.28 0.6741
6.16 0.7128
8.86 0.7528
12.74 0.7979
18.33 0.8398
26.37 0.8792
37.93 0.9141
54.56 0.9426
78.48 0.9674
112.88 0.9877
162.38 1.0048
233.57 1.0188
335.98 1.0309
483.29 1.0423
695.19 1.0601
1000 1.0716
};
\addlegendentry{A-ADMM ($\bar n = {1}$)}

% --- A-ADMM (n_OZF=4) ---
\addplot [ozf4, thick, mark=*, dashed, mark size=1.25pt, mark options={solid, fill=ozf4, draw=ozf4}]
table {%
1 0.5001953125
1.43844988828766 0.545263671875
2.06913808111479 0.589697265625
2.97635144163132 0.63349609375
4.28133239871939 0.67412109375
6.15848211066026 0.712841796875
8.85866790410082 0.75283203125
12.7427498570313 0.7947265625
18.3298071083244 0.83154296875
26.3665089873036 0.861376953125
37.9269019073225 0.8861328125
54.5559478116852 0.907080078125
78.4759970351461 0.923583984375
112.883789168469 0.93818359375
162.377673918872 0.950244140625
233.572146909012 0.960400390625
335.981828628378 0.96865234375
483.293023857175 0.975634765625
695.192796177561 0.98134765625
1000 0.98642578125
};
\addlegendentry{A-ADMM ($\bar n = {4}$)}

% --- Optimal rho (solid, same red) ---
\addplot [refgray, line width=1pt]
table {%
1.000000       0.000000
1.63789371     0.28099203
2.6826958      0.47029191
4.39397056     0.60405492
7.19685673     0.70105667
11.78768635    0.77260740
19.30697729    0.82603977
31.6227766     0.86631946
51.79474679    0.89690938
84.83428982    0.92027735
138.94954944   0.93821226
227.58459261   0.95202889
372.75937203   0.96270468
610.54022966   0.97097320
1000.000000    0.97738929
};
\addlegendentry{Opt. $\rho$ of NM method.}

\end{axis}
\end{tikzpicture}
    \caption{\rev{Comparison of $\rho$ values and influence of O'Shea-Zames-Falb order for A-ADMM with NM tuning.}}
    \label{fig:NM_study}
\end{figure}

\vspace{-2mm}
%%%%%%%%%%%%%%%%%%%%%%%%%%%%%%%%%%%%%%%%%%%%%%%%%%%%%

%%%%%%%%%%%%%%%%%%%%%%%%%%%%%%%%%%%%%%%%%%%%%%%%%%%%%
\subsection{Grid Search Results}
\label{sec:grid search}
\vspace{-3mm}

\meisam{Next, we perform a grid search to identify optimal 
$(\nu_1, \nu_2)$ configurations. 
For each $\kappa$, we set up a fixed grid from which
$(\nu_1, \nu_2)$ is selected as the pair minimizing $\rho$ via 
bisection on Theorem~\ref{thm:worst-case cng thm}. The procedure is repeated for the over-relaxed version (OR-A-ADMM) with an additional grid over 
$\alpha \in (0, 4]$.}

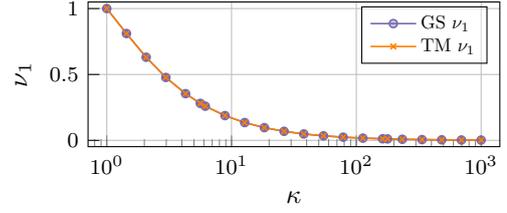
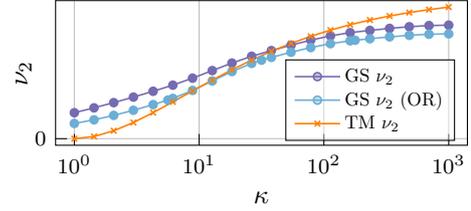
\begin{figure} 
    \centering
    \begin{subfigure}{0.97\columnwidth}
        \centering
        \begin{tikzpicture}

\definecolor{skyblue}{RGB}{107,174,214}
\definecolor{purple}{RGB}{117,107,177}
\definecolor{orange}{RGB}{255,127,0}

\begin{axis}[
  width=7cm,
    height=3.5cm,
  every axis/.append style={},     % <--- global font size
  tick label style={font=\footnotesize},      % <--- smaller tick labels
  legend cell align={left},
  legend style={
  fill opacity=0.4,
  draw opacity=1,
  text opacity=1,
  nodes={scale=0.7, transform shape},
  at={(0.97,0.97)},
  anchor=north east,
  draw=black
},
log basis x={10},
tick align=inside,
tick pos=left,
% x grid style={darkgray176},
xlabel={\(\displaystyle \kappa\)},
xlabel style={
  at={(axis description cs:0.5,-0.04)},
  anchor=north
},
xmajorgrids,
xmin=0.707945784384138, xmax=1412.53754462275,
xmode=log,
% y grid style={darkgray176},
ylabel={$\nu_1$},
ymajorgrids,
ymin=-0.0479332039154318, ymax=1.04990158113883,
% ylabel style={yshift=-4mm},           % <- put the y-label closer to the axis
ytick style={color=black}
]
\addplot [semithick, purple, mark=*, mark size=1.5, mark options={solid}]
table {%
1.00 1.0000
1.44 0.8446
2.07 0.7101
2.98 0.5951
4.28 0.4974
6.16 0.4149
8.86 0.3456
12.74 0.2874
18.33 0.2387
26.37 0.1981
37.93 0.1643
54.56 0.1363
78.48 0.1131
112.88 0.0941
162.38 0.0788
233.57 0.0668
335.98 0.0578
483.29 0.0519
695.19 0.0491
1000.00 0.0500
};
\addlegendentry{GS $\nu_1$}

\addplot [semithick, skyblue, mark=*, mark size=1.5, mark options={solid}]
table {%
1.00 1.0000
1.44 0.8197
2.07 0.6757
2.98 0.5590
4.28 0.4636
6.16 0.3851
8.86 0.3203
12.74 0.2666
18.33 0.2220
26.37 0.1849
37.93 0.1541
54.56 0.1284
78.48 0.1071
112.88 0.0892
162.38 0.0744
233.57 0.0620
335.98 0.0517
483.29 0.0431
695.19 0.0359
1000.00 0.0300
};
\addlegendentry{{GS $\nu_1$ (OR)}}

%\addplot [semithick, orange, mark=x, mark size=1.5, mark options={solid}]
%table {%
%1 1
%1.43844988828766 0.810746197016624
%2.06913808111479 0.630604219085972
%2.97635144163132 0.477215353265881
%4.28133239871939 0.354260504649556
%5.62341325190349 0.280666461074539
%6.15848211066026 0.259323456540615
%8.85866790410082 0.187840676429615
%12.7427498570313 0.134968067581669
%18.3298071083244 0.0963691457663391
%26.3665089873036 0.0684676039925656
%% 31.6227766016838 0.0576221399514641
%37.9269019073225 0.0484516855758878
%54.5559478116852 0.0341779852938119
%78.4759970351461 0.024047049825775
%112.883789168469 0.0168835535847299
%162.377673918872 0.0118336711974634
%177.827941003892 0.0108251300003784
%233.572146909012 0.0082825291213189
%335.981828628378 0.00579032520934377
%483.293023857175 0.00404415566550277
%695.192796177561 0.00282234382876364
%1000 0.00196837722339832
%};
%\addlegendentry{TM $\nu_1$}

\end{axis}

\end{tikzpicture} 
        \vspace{-1mm}
        \caption{Results for $\nu_1$.}
        \label{fig:grid_vs_tmm}
    \end{subfigure}
    \hfill
    \begin{subfigure}{0.97\columnwidth}
        \centering
        \begin{tikzpicture}

\definecolor{skyblue}{RGB}{107,174,214}
\definecolor{purple}{RGB}{117,107,177}
\definecolor{orange}{RGB}{255,127,0}

\begin{axis}[
  width=7cm,
    height=3.5cm,
  every axis/.append style={},     % <--- global font size
  tick label style={font=\footnotesize},      % <--- smaller tick labels
  legend cell align={left},
  legend style={
  fill opacity=0.4,
  draw opacity=1,
  text opacity=1,
  nodes={scale=0.7, transform shape},
  at={(0.03,0.97)},
  anchor=north west,
  draw=black
},
log basis x={10},
tick align=inside,
tick pos=left,
% x grid style={darkgray176},
xlabel={\(\displaystyle \kappa\)},
xlabel style={
  at={(axis description cs:0.5,-0.04)},
  anchor=north
},
xmajorgrids,
xmin=0.707945784384138, xmax=1412.53754462275,
xmode=log,
% y grid style={darkgray176},
ylabel={$\nu_2$},
ymajorgrids,
ymin=-0.0479332039154318, ymax=1.04990158113883,
% ylabel style={yshift=-4mm},           % <- put the y-label closer to the axis
ytick style={color=black}
]
\addplot [semithick, purple, mark=*, mark size=1.5, mark options={solid}]
table {%
1.00 0.1153
1.44 0.1258
2.07 0.1300
2.98 0.1649
4.28 0.2160
6.16 0.2738
8.86 0.3325
12.74 0.3889
18.33 0.4410
26.37 0.4880
37.93 0.5297
54.56 0.5661
78.48 0.5977
112.88 0.6249
162.38 0.6481
233.57 0.6679
335.98 0.6847
483.29 0.6988
695.19 0.7108
1000.00 0.7209
};
\addlegendentry{{GS $\nu_2$}}

\addplot [semithick, skyblue, mark=*, mark size=1.5, mark options={solid}]
table {%
1.00 0.0020
1.44 0.0106
2.07 0.0164
2.98 0.0306
4.28 0.0606
6.16 0.0978
8.86 0.1370
12.74 0.1753
18.33 0.2108
26.37 0.2429
37.93 0.2713
54.56 0.2961
78.48 0.3175
112.88 0.3358
162.38 0.3514
233.57 0.3647
335.98 0.3759
483.29 0.3853
695.19 0.3933
1000.00 0.4000
};
\addlegendentry{{GS $\nu_2$ (OR)}}

%\addplot [semithick, orange, mark=x, mark size=1.5, mark options={solid}]
%table {%
%1 0
%1.43844988828766 0.0150663196972622
%2.06913808111479 0.0548064100505549
%2.97635144163132 0.111862896350163
%4.28133239871939 0.17999552003581
%6.15848211066026 0.254073617396337
%8.85866790410082 0.330034527763297
%12.7427498570313 0.40480447069959
%18.3298071083244 0.476187514014118
%26.3665089873036 0.542733805739324
%37.9269019073225 0.603599997567166
%54.5559478116852 0.658413531423883
%78.4759970351461 0.707149460868785
%112.883789168469 0.75002495328939
%162.377673918872 0.787413433748117
%233.572146909012 0.819777929437706
%335.981828628378 0.847621653145216
%483.293023857175 0.871453112156774
%695.192796177561 0.891762844158643
%1000 0.909009056474822
%};
%\addlegendentry{TM $\nu_2$}1

\end{axis}

\end{tikzpicture} 
        % \vspace{-2mm} 
        \caption{Results for $\nu_2$. }
    \end{subfigure}
    \vspace{-3mm}
    \caption{Results of the grid search for ($\nu_1$, $\nu_2$).}
    \label{fig:GS-results}
\end{figure}
Fig.~\ref{fig:GS-results} shows the results for $\nu_{1}$ and $\nu_{2}$.
\rev{Applying symbolic regression \citep{pysr} reveals that $\nu_{1}$ and $\nu_{2}$ can be approximated
by}
\meisam{\begin{align*}
\nu_1 &\approx \frac{2.26}{\log(\kappa + 5.49)}, &\quad \nu_2 &\approx 0.27\sqrt{\log(\kappa)}, \\
\nu_1^{\mathrm{OR}} &= \nu_1, &\quad \nu_2^{\mathrm{OR}} &\approx 0.079\,\sqrt[4]{\kappa},
\end{align*}}for the nominal and over-relaxed settings, respectively.
The over-relaxation parameter $\alpha$ was found to be optimal for $\alpha \approx 1.45$ across all condition numbers. 

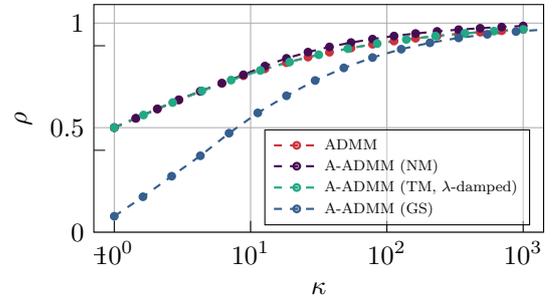
\begin{figure}
  \centering
  \begin{tikzpicture}

\definecolor{darkgray176}{RGB}{176,176,176}
\definecolor{goldenrod18818934}{RGB}{188,189,34}
\definecolor{green}{RGB}{0,128,0}
\definecolor{lightgray204}{RGB}{204,204,204}
\definecolor{purple}{RGB}{128,0,128}
\definecolor{darkgray176}{RGB}{176,176,176}
\definecolor{mycolor1}{RGB}{204,51,63}   % red
\definecolor{mycolor2}{RGB}{68,12,84}    % purple
\definecolor{mycolor3}{RGB}{53,95,141}   % blue
\definecolor{mycolor4}{RGB}{34,168,132}  % teal
\definecolor{mycolor5}{RGB}{230,142,36}
\definecolor{ozf4}{RGB}{8,48,107}

\begin{axis}[
    width=7cm,
    height= 5cm,
    % every axis/.append style={font=\small},
    % tick label style={font=\small},
    legend style={
      fill opacity=0.4,
      draw opacity=1,
      text opacity=1,
      font=\small,
      nodes={scale=0.7, transform shape},
      at={(0.97,0.03)},
      anchor=south east,
      legend cell align={left}
    },
    log basis x={10},
    tick align=inside,
    tick pos=left,
    x grid style={darkgray176},
    xlabel={$\kappa$},
    xlabel style={
      at={(axis description cs:0.5,-0.04)},
      anchor=north
    },
    xmajorgrids,
    xmin=0.707945784384138, xmax=1412.53754462275,
    xmode=log,
    xtick style={color=black},
    xtick={0.01,0.1,1,10,100,1000,10000,100000,1000000},
    xticklabels={
      \(\displaystyle {10^{-2}}\),
      \(\displaystyle {10^{-1}}\),
      \(\displaystyle {10^{0}}\),
      \(\displaystyle {10^{1}}\),
      \(\displaystyle {10^{2}}\),
      \(\displaystyle {10^{3}}\),
      \(\displaystyle {10^{4}}\),
      \(\displaystyle {10^{5}}\),
      \(\displaystyle {10^{6}}\)
    },
    y grid style={darkgray176},
    ylabel={$\rho$},
    ymajorgrids,
    ymin=0.25, ymax=1.08875,
    ytick style={color=black, yshift=-3mm}
]
% \addplot [mycolor1, thick, mark=*, dashed, mark size=1.25pt, mark options={solid, fill=mycolor1, draw=mycolor1}]
% table {%
% 1.00 0.5002
% 1.44 0.5440
% 2.07 0.5891
% 2.98 0.6322
% 4.28 0.6735
% 6.16 0.7122
% 8.86 0.7484
% 12.74 0.7814
% 18.33 0.8106
% 26.37 0.8373
% 37.93 0.8607
% 54.56 0.8811
% 78.48 0.8988
% 112.88 0.9141
% 162.38 0.9274
% 233.57 0.9388
% 335.98 0.9483
% 483.29 0.9566
% 695.19 0.9636
% 1000.00 0.9699
% };
% \addlegendentry{ADMM}

\addplot [mycolor1, thick, mark=*, dashed, mark size=1.25pt, mark options={solid, fill=mycolor1, draw=mycolor1}]
table {%
1 0.5002
1.47 0.5465
2.15 0.5941
3.16 0.6398
4.64 0.6817
6.81 0.7230
10 0.7598
14.68 0.7935
21.54 0.8227
31.62 0.8493
46.42 0.8722
68.13 0.8925
100 0.9096
146.78 0.9242
215.44 0.9363
316.23 0.9471
464.16 0.9560
681.29 0.9636
1000 0.9699
};
\addlegendentry{ADMM}

\addplot [mycolor4, thick, mark=*, dashed, mark size=1.25pt, mark options={solid, fill=mycolor4, draw=mycolor4}]
table {%
1 0.5002
1.44 0.5440
2.07 0.5884
2.98 0.6322
4.28 0.6735
6.16 0.7128
8.86 0.7484
12.74 0.7814
18.33 0.8106
26.37 0.8373
37.93 0.8607
54.56 0.8906
78.48 0.9217
112.88 0.9502
162.38 0.9750
233.57 0.9959
335.98 1.0137
483.29 1.0296
695.19 1.0480
1000 1.0645
};
\addlegendentry{A-ADMM (TM)}

\addplot [ozf4, thick, mark=*, dashed, mark size=1.25pt, mark options={solid, fill=ozf4, draw=ozf4}]
table {%
1 0.5001953125
1.43844988828766 0.545263671875
2.06913808111479 0.589697265625
2.97635144163132 0.63349609375
4.28133239871939 0.67412109375
6.15848211066026 0.712841796875
8.85866790410082 0.75283203125
12.7427498570313 0.7947265625
18.3298071083244 0.83154296875
26.3665089873036 0.861376953125
37.9269019073225 0.8861328125
54.5559478116852 0.907080078125
78.4759970351461 0.923583984375
112.883789168469 0.93818359375
162.377673918872 0.950244140625
233.572146909012 0.960400390625
335.981828628378 0.96865234375
483.293023857175 0.975634765625
695.192796177561 0.98134765625
1000 0.98642578125
};
\addlegendentry{A-ADMM (NM)}

\addplot [mycolor5, thick, mark=*, dashed, mark size=1.25pt, mark options={solid, fill=mycolor5, draw=mycolor5}]
table {%
1.00 0.2977
1.44 0.3578
2.07 0.4180
2.98 0.4773
4.28 0.5345
6.16 0.5886
8.86 0.6391
12.74 0.6856
18.33 0.7277
26.37 0.7656
37.93 0.7992
54.56 0.8288
78.48 0.8548
112.88 0.8775
162.38 0.8973
233.57 0.9146
335.98 0.9299
483.29 0.9401
695.19 0.9504
1000.00 0.9573
};
\addlegendentry{A-ADMM (GS)}

\end{axis}

\end{tikzpicture}
  \vspace{-2mm}
\caption{Comparison of optimal $\rho$ values for ADMM and A-ADMM with different parameter configurations.}
\label{fig:A-ADMM}
\end{figure}

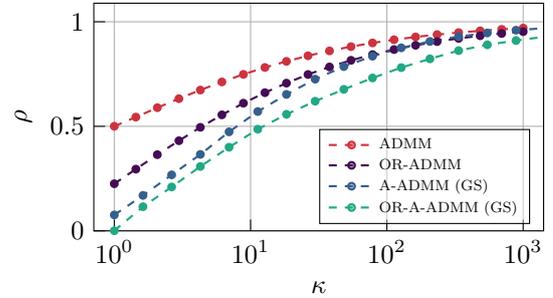
\begin{figure}
    \centering
    % This file was created with tikzplotlib v0.10.1.
\begin{tikzpicture} 

\definecolor{darkgray176}{RGB}{176,176,176}
\definecolor{goldenrod18818934}{RGB}{188,189,34}
\definecolor{green}{RGB}{0,128,0}
\definecolor{lightgray204}{RGB}{204,204,204}
\definecolor{purple}{RGB}{128,0,128}
\definecolor{darkgray176}{RGB}{176,176,176}
\definecolor{mycolor1}{RGB}{204,51,63}   % red
\definecolor{mycolor2}{RGB}{68,12,84}    % purple
\definecolor{mycolor3}{RGB}{53,95,141}   % blue
\definecolor{mycolor4}{RGB}{34,168,132}  % teal
\definecolor{mycolor5}{RGB}{230,142,36}

\begin{axis}[
    width=7cm,
    height=5cm,
    % every axis/.append style={font=\small},
    % tick label style={font=\small},
    legend style={
      fill opacity=0.4,
      draw opacity=1,
      text opacity=1,
      font=\small,
      nodes={scale=0.7, transform shape},
      at={(0.97,0.03)},
      anchor=south east,
      legend cell align={left}
    },
    log basis x={10},
    tick align=inside,
    tick pos=left,
    x grid style={darkgray176},
    xlabel={$\kappa$},
    xlabel style={
      at={(axis description cs:0.5,-0.04)},
      anchor=north
    },
    xmajorgrids,
    xmin=0.707945784384138, xmax=1412.53754462275,
    xmode=log,
    xtick style={color=black},
    xtick={0.01,0.1,1,10,100,1000,10000,100000,1000000},
    xticklabels={
      \(\displaystyle {10^{-2}}\),
      \(\displaystyle {10^{-1}}\),
      \(\displaystyle {10^{0}}\),
      \(\displaystyle {10^{1}}\),
      \(\displaystyle {10^{2}}\),
      \(\displaystyle {10^{3}}\),
      \(\displaystyle {10^{4}}\),
      \(\displaystyle {10^{5}}\),
      \(\displaystyle {10^{6}}\)
    },
    y grid style={darkgray176},
    ylabel={$\rho$},
    ymajorgrids,
    ymin=0.15, ymax=1.0,
    ytick style={color=black}
]
\addplot [mycolor1, thick, mark=*, dashed, mark size=1.25pt, mark options={solid, fill=mycolor1, draw=mycolor1}]
table {%
1.00 0.5002
1.44 0.5440
2.07 0.5891
2.98 0.6322
4.28 0.6735
6.16 0.7122
8.86 0.7484
12.74 0.7814
18.33 0.8106
26.37 0.8373
37.93 0.8607
54.56 0.8811
78.48 0.8988
112.88 0.9141
162.38 0.9274
233.57 0.9388
335.98 0.9483
483.29 0.9566
695.19 0.9636
1000.00 0.9699
};
\addlegendentry{ADMM}

\addplot [mycolor2, thick, mark=*, dashed, mark size=1.25pt, mark options={solid, fill=mycolor2, draw=mycolor2}]
table {%
1.00 0.3009
1.44 0.3637
2.07 0.4259
2.98 0.4862
4.28 0.5440
6.16 0.5979
8.86 0.6481
12.74 0.6938
18.33 0.7351
26.37 0.7719
37.93 0.8049
54.56 0.8334
78.48 0.8582
112.88 0.8798
162.38 0.8982
233.57 0.9141
335.98 0.9280
483.29 0.9395
695.19 0.9490
1000.00 0.9572
};
\addlegendentry{OR-ADMM}
\addplot [mycolor5, thick, mark=*, dashed, mark size=1.25pt, mark options={solid, fill=mycolor5, draw=mycolor5}]
table {%
1.00 0.2977
1.44 0.3578
2.07 0.4180
2.98 0.4773
4.28 0.5345
6.16 0.5886
8.86 0.6391
12.74 0.6856
18.33 0.7277
26.37 0.7656
37.93 0.7992
54.56 0.8288
78.48 0.8548
112.88 0.8775
162.38 0.8973
233.57 0.9146
335.98 0.9299
483.29 0.9401
695.19 0.9504
1000.00 0.9573
};
\addlegendentry{A-ADMM (GS)}
\addplot [mycolor4, thick, mark=*, dashed, mark size=1.25pt, mark options={solid, fill=mycolor4, draw=mycolor4}]
table {%
1.00 0.1759
1.44 0.2366
2.07 0.3004
2.98 0.3668
4.28 0.4342
6.16 0.5004
8.86 0.5634
12.74 0.6215
18.33 0.6738
26.37 0.7199
37.93 0.7599
54.56 0.7944
78.48 0.8239
112.88 0.8490
162.38 0.8706
233.57 0.8891
335.98 0.9054
483.29 0.9201
695.19 0.9340
1000.00 0.9477
};
\addlegendentry{OR-A-ADMM (GS)}
\end{axis}

\end{tikzpicture}
     \vspace{-2mm}
    \caption{Comparison of optimal $\rho$ values for nominal and over-relaxed ADMM versions.}
    \label{fig:OR-A-ADMM}
\end{figure}

Fig. \ref{fig:A-ADMM} displays the resulting convergence rates obtained from the grid search, in the context of the vanilla and accelerated schemes. By design, A-ADMM with GS parameters gives the best linear convergence rate among all settings. Fig. \ref{fig:OR-A-ADMM} complements the results with the over-relaxed variants. 
\meisam{We observe that the over-relaxed versions consistently improve upon their non-over-relaxed counterparts across the full range of condition numbers. Combining over-relaxation with the GS-tuned acceleration yields the fastest certified linear convergence rate among all configurations considered in this work.}

\vspace{-2mm}
\section{Case Study}
\vspace{-3mm}
We test the A-ADMM schemes designed based on the theoretical analyses in the previous section on the $\ell_1$-regularized
least-squares problem, also known as {LASSO regression}, and compare them with available benchmark algorithms \citep{boydAdmm}.
We consider the problem
\begin{align} 
    \min_{x, z \in \mathbb{R}^{500}} \quad & \tfrac{1}{2} \|F x - b\|_2^2 + \tau \|z\|_1 \nonumber \\
    \text{s.t.} \quad & x = z,
    \label{eq:lasso_admm}
\end{align}
\meisam{where $F \in \mathbb{R}^{600 \times 500}$ is chosen as a full rank 
%random matrix with normalized columns.
matrix with entries that are first sampled from an isotropic Gaussian distribution $\mathcal{N}(0,1)$, whose columns are then normalized. 
The vector $b$ is generated as $b = F w^0 + \varepsilon$, where $w^0 \in \mathbb{R}^{500}$ is chosen sparse, containing 250
non-zeros drawn from $\mathcal{N}(0,1)$ 
and $\varepsilon \sim \mathcal{N}(0,10^{-3}I_{600})$. The regularization parameter
$\tau > 0$ is set to $\tau=1$.} 
Problem \eqref{eq:lasso_admm} fits 
%the general problem 
\eqref{constrained optimization}, where the first term is $m$-strongly convex and $L$-smooth with 
$m = \underline{\sigma}^2(F)$ and $L = \bar{\sigma}^2(F)$.

Applying the accelerated over-relaxed ADMM results in
\begin{subequations}\label{Lasso-Proposed-A-ADMM}
\begin{align}
    x_{k+1} &= (\nu_1 F^\top F + I)^{-1} \left( -\hat{z}_k - \hat{\lambda}_k + \nu_1 F^\top b \right) \label{eq:oradmm_x}, \\
    z_{k+1} &= \mathcal{S}_{\tau \nu_1 } \left( - \alpha\, x_{k+1} - (\alpha - 1)\hat{z}_k - \hat{\lambda}_k \right) \label{eq:oradmm_z}, \\
    \lambda_{k+1} &= \alpha\, x_{k+1} + (1 - \alpha)\,\hat{z}_k - z_{k+1} + \hat{\lambda}_k \label{eq:oradmm_u}, 
\end{align}
\end{subequations}
with $\hat z_k$ and $\hat \lambda_k$ as in \eqref{eq:ADMM_accel}, and $\mathcal{S}_{\tau \nu_1}$ being the {soft-thresholding} operator,  which is (element-wise) defined as
\begin{align}\label{eq:soft_thresholding}
\bigl[\mathcal{S}_{\tau \nu_1}(\mathbf{y})\bigr]_i
=
\operatorname{sign}(\mathbf{y}_i)\,
\max\bigl(|\mathbf{y}_i|-\tau \nu_1,0\bigr).
\end{align}
The convergence is assessed through the normalized iterate error
\begin{align}
    \Delta_k := 
    \frac{\|x_k - x^\star\|_2}{\|x_0 - x^\star\|_2},
    \label{eq:iterate_subopt}
\end{align}
where $x^\star$ is the solution of~\eqref{eq:lasso_admm}. 

\begin{figure}
    \centering
    \input{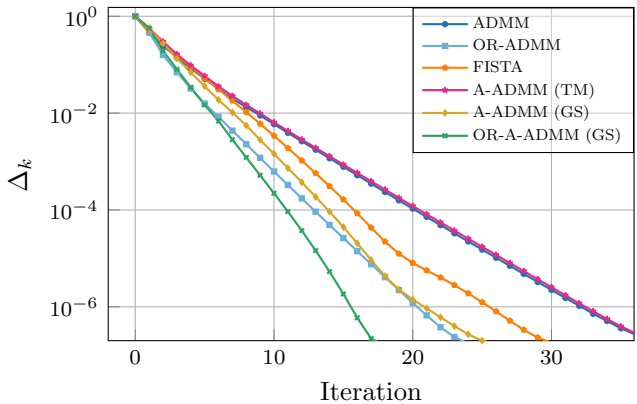} 
    \vspace{-2mm}
    \caption{Empirical performance of different A-ADMM schemes and other benchmark algorithms applied to a random instance of LASSO regression.}
    \label{fig:final_summary}
\end{figure}

\meisam{Specifically, we compare A-ADMM using NM and GS parameter selections against the widely used LASSO benchmark solver FISTA~\citep{FISTA} and the non-accelerated baselines ADMM and OR-ADMM.} The iterate error computed on one random instance of \eqref{eq:lasso_admm} is shown in Fig.~\ref{fig:final_summary}, while Table~\ref{tb:average_iters} summarizes the average number of iterations required to reach \meisam{$\Delta_k < 10^{-4}$} over 200 random instances.
\vspace{-2mm}

\begin{table}
\centering
\caption{\meisam{Average iterations to reach $\Delta_k < 10^{-4}$ over 200 random LASSO instances.}}
\label{tb:average_iters}
\begin{tabular}{lc}
\toprule
Method & Average iterations \\ 
\midrule
ADMM & 197.66 \\
OR-ADMM & 163.63 \\
A-ADMM (NM) & 137.30 \\
FISTA & 136.71 \\
A-ADMM (GS) & 108.03 \\
OR-A-ADMM (GS) & 98.77 \\
\bottomrule
\end{tabular}
\end{table}

Both Fig.~\ref{fig:final_summary} and Table~\ref{tb:average_iters} show that the performance of the A-ADMM methods differs profoundly based on their parameter selection. 
We also observe that the accelerated versions improve the convergence speed substantially over the vanilla \meisam{baselines}. Finally, OR-A-ADMM with GS tuning consistently outperforms all other methods. 
Notably, the ordering of algorithms evaluated on the empirical performance and evaluated on the worst-case convergence rates (cf. Fig.~\ref{fig:OR-A-ADMM}) stays the same, suggesting that the worst-case convergence rate is a suitable metric for algorithm comparison.
\vspace{-3mm}
\section{Conclusions}
\vspace{-3mm}

\meisam{This paper extended the IQC-based algorithm analysis framework to certify worst-case linear convergence rates of A-ADMM for any choice of the algorithm parameters, thereby enabling a systematic investigation of A\nobreakdash-ADMM performance under different parameter selections. We showed that dynamic O'Shea-Zames-Falb multipliers are essential to certify convergence in the accelerated setting.} We further proposed different parameter tunings and showed that the additional momentum parameters allow \meisam{A-ADMM} to outperform vanilla ADMM. %The fastest version,
The fastest configuration, in terms of both certified convergence rate and empirical performance, was found using a grid search over the parameter space, giving the fastest version of A-ADMM in the strongly-convex-smooth setting. %We moreover showed the close connection between A-ADMM parameters and accelerated gradient methods from unconstrained optimization.
Future work could focus on closed-form characterizations of the optimal A-ADMM parameters and on a systematic A\nobreakdash-ADMM synthesis procedure.

\vspace{-1mm}

\bibliography{references}
% \bibliography{main}
\end{document}